\documentclass[10pt, a4paper, reqno, twoside]{article}

\usepackage{amssymb}
\usepackage{amsmath} 
\usepackage{amsthm}

\usepackage{rotating}
\usepackage{graphicx}

\newtheorem{thm}{Theorem}[section]

\newtheorem{lemma}[thm]{Lemma}

\newtheorem{Remark}[thm]{Remark}

\newtheorem{cor}[thm]{Corollary}

\def\medint{-\kern  -,375cm\int}

\def\[{\bigl(}  \def\]{\bigl)}

\def\rz{\ifmmode{I\hskip -3.2pt R}
        \else{\hbox{$I\hskip -3.2pt R$}}\fi} 

\begin{document}

\title{ On the geometry of the $p$-Laplacian operator}

\author{B.~Kawohl \&  J.~Hor\'ak }

\maketitle

\bigskip\bigskip
\begin{abstract} The $p$--Laplacian operator $\Delta_pu={\rm div }\left(|\nabla u|^{p-2}\nabla u\right)$  is not uniformly elliptic for any $p\in(1,2)\cup(2,\infty)$ and degenerates even more when $p\to \infty$ or $p\to 1$. In those two cases the Dirichlet and eigenvalue problems associated with
the $p$-Laplacian
lead to intriguing geometric questions, because their limits for $p\to\infty$ or $p\to 1$ can be characterized by the geometry of $\Omega$. In this little survey I recall some well-known results on eigenfunctions of the classical 2-Laplacian and elaborate on their extensions to general $p\in[1,\infty]$.
 We report also on results concerning the normalized or game-theoretic $p$--Laplacian 
$$\Delta_p^Nu:=\tfrac{1}{p}|\nabla u|^{2-p}\Delta_pu=\tfrac{1}{p}\Delta_1^Nu+\tfrac{p-1}{p}\Delta_\infty^Nu$$
and its parabolic counterpart $u_t-\Delta_p^N u=0$. These equations are homogeneous of degree 1 and $\Delta_p^N$ is uniformly elliptic for any $p\in  (1,\infty)$. In this respect it is more benign than the $p$-Laplacian, but it is not of divergence type.

\end{abstract}

\section{Introduction}

In intrinsic coordinates, orthogonal and tangential to level surfaces of a function $u$, one can write the linear Laplacian as
$$\Delta u= u_{x_1x_1}+\ldots+ u_{x_nx_n}=u_{\nu\nu}+u_\nu\ {\rm div}(\nu)$$
where $\nu(x)=-\tfrac{\nabla u(x)}{|\nabla u(x)|}$ is direction of steepest descent.
In fact, 
$${\rm div}(\nu)=-\frac{\Delta u}{|\nabla u|}+\frac{u_{x_i}u_{x_j}u_{x_i x_j}}{|\nabla u|^3}=-\frac{\Delta u}{|\nabla u|}+\frac{u_{\nu\nu}}{|\nabla u|}$$ 
so that $\Delta u=u_{\nu\nu}-|\nabla u|\ {\rm div}(\nu)=u_{\nu\nu}+u_\nu\ {\rm div}(\nu)$ or 
\begin{equation}\label{lap}
\Delta u=u_{\nu\nu}+(n-1)\,H\, u_\nu
\end{equation}
with $H$ denoting mean curvature of a level set  of $u$.  This observation is more or less familiar to anyone who ever calculated the Laplacian in polar coordinates. If $u(x)=v(|x|)=v(r)$ is radial and  radially decreasing, then $\quad\Delta v=v_{rr}+\tfrac{n-1}{r}v_r$ and $\tfrac{1}{r}$ is the mean curvature of the sphere of radius $r$.
In a similar way for $p\in[1,\infty)$ one can write the {\bf $p$-Laplace operator} as
\begin{eqnarray}\label{plap}
\Delta_p u&=&{\rm div}\left(|\nabla u|^{p-2}\nabla u\right)=|\nabla u|^{p-2}[\Delta u+(p-2)u_{\nu\nu}]\cr
&=&|\nabla u|^{p-2}\left[(p-1)u_{\nu\nu}+(n-1)Hu_\nu\right]\end{eqnarray}
and the {\bf normalized} or {\bf game-theoretic $p$-Laplace operator} as
\begin{eqnarray}\label{nplap}
\Delta_p^N u&=&\tfrac{1}{p}|\nabla u|^{2-p}{\rm div}\left(|\nabla u|^{p-2}\nabla u\right)\cr
&=&\tfrac{p-1}{p}u_{\nu\nu}+\tfrac{1}{p}(n-1)Hu_\nu=\tfrac{p-1}{p}\Delta_\infty^Nu+\tfrac{1}{p}\Delta_1^Nu\ .
\end{eqnarray}
Observe $\quad\Delta_\infty^Nu=u_{\nu\nu}$, $\quad \Delta_2^Nu=\tfrac{1}{2}\Delta u\quad$ and $\quad\Delta_1^Nu=|\nabla u|{\rm div}(\tfrac{\nabla u}{|\nabla u|})$. 

\section{Dirichlet problems for $-\Delta_pu=0$ and $-\Delta_pu=1$}\label{Dirichletsection}\setcounter{equation}{0}

For $p\in(1,\infty)$ these problems are well understood. Unique solutions can be found by variational methods involving the strictly convex functionals 
$$\int_\Omega\tfrac{1}{p}|\nabla v|^p\ dx\qquad\hbox{ or }\qquad\int_\Omega\left\{\tfrac{1}{p}|\nabla v|^p-v\right\}\ dx\qquad\hbox{ on }\ W^{1,p}_0(\Omega)+g.$$

Let us recall that the limit $p\to\infty$ of $p$-harmonic functions $-\Delta_pu=0$ in $\Omega$, $u=g$ on $\partial\Omega$
solves 
\begin{equation}\label{dirinfty}
-\Delta_\infty u=-\sum_{i,j=1}^nu_{x_i}u_{{x_i}{x_j}}u_{x_j}=0\quad\hbox{ in }\Omega, \qquad\qquad u=g\quad \hbox{ on }\partial\Omega
\end{equation}
in the sense of viscosity solutions. Most of us had the privilege of playing with dry sand (and not virtual toys) as toddlers, and since then we are familiar with $\infty$-harmonic functions. This equation models the shape of sandpiles. When the slope of a sandpile becomes to large, the sand starts rolling down. This observation caused Gunnar Aronsson to derive the corresponding equation in 
\cite{Ar}. A famous explicit solution in two dimensions is $u(x,y)=x^{4/3}-y^{4/3}$. The saddle that it describes in the origin can be visualized as the place where two initially separate and conical, but growing sandpiles meet \cite{KS}. This solution is not of class $C^2$ along the coordinate axes, and so one has to interpret the differential equation for $p=\infty$ \lq\lq in the viscosity sense\rq\rq$ $ \cite{CIL, K}. It is remarkable, that solutions to (\ref{dirinfty}) are unique \cite{Je}, although the equation is not strongly elliptic, while the minimizers of the limiting functional can be non-unique. 

Incidentally, (viscosity) solutions to $-\Delta_\infty=f$ in $\Omega$, $u=0$ on $\partial\Omega$,
continue to exist  and are still unique, if $f$ has only one sign, but there are counterexamples to uniqueness for $f$ changing sign, see \cite{LW0}

The limit $p\to1$ of $p$-harmonic functions as $p\to 1$ was investigated by Juutinen in \cite{J1}. Formally, from (\ref{plap}), we expect 1-harmonic functions to be solutions of 
\begin{equation}\label{dir1}
-\Delta_1 u=-{\rm div}\left(\frac{\nabla u}{|\nabla u|}\right)=-(n-1)H=0\ \hbox{ in }\Omega\qquad\qquad u=g\hbox{ on }\ \partial\Omega.
\end{equation}
An example of nonuniqueness for solutions to (\ref{dir1}) is depicted in \cite{K}. However, as pointed out in \cite{J1}, the limiting functional has a unique minimizer.

The nonlinear torsion problem
\begin{equation}\label{ptorsion}
-\Delta_pu=1\quad\hbox{ in }\Omega, \qquad\qquad u=0\quad\hbox{ on }\partial\Omega
\end{equation} 
and its limits $p\to\infty$ and $p\to 1$ were investigated in \cite{K0, BBM}. In contrast to the case of vanishing right hand side, the limiting equation for $p\to\infty$ is not (~as one might expect from (\ref{dirinfty})~) $-\Delta_\infty u=1$, but instead $|\nabla u|=1$.

For the sake of exposition l want to give the proof of  the slightly weaker result, that the limiting equation is characterized by
\begin{equation}\label{torsioninfty}
\min\{|\nabla u|-1, -\Delta_\infty u\}=0 \quad\hbox{ in }\Omega,\qquad u=0\quad\hbox{ on }\partial\Omega.
\end{equation}
Let us see why. From the weak formulation of $-\Delta_p u_p=0$ in $\Omega$, $u_p=0$ on $\partial\Omega$ and Poincar\'e's inequality one may conclude that for a fixed $s>n$ and any $p>s$
$$\lim_{p\to\infty}||\nabla u_p||_s\to |\Omega|^{1/s}.$$ Therefore the family $\{ u_p\}_{p>s}$ is uniformly bounded in $W^{1,s}_0(\Omega)$ and after passing to a subsequence it converges uniformly to a H\"older-continuous limit $u_\infty$. Now suppose that a $C^2$ test function $v$ touches $u_\infty$ from above in $x_\infty\in\Omega$. Then also $\tilde{v}=v+M|x-x_\infty|^4$ touches there (for large positive $M$). Moreover, $\tilde{v}-u_p$ has a local minimum in $x_p$ near $x_\infty$, and thus $\hat{v}(x):=\tilde{v}(x)+u_p(x_p)-\tilde{v}(x_p)$ touches $u_p$ in $x_p$ from above. Now $u_p$ is a viscosity subsolution of (\ref{plap}), so
\begin{equation}
-|\nabla \hat{v}|^{p-4}\left(|\nabla\hat{v}|^2\Delta\hat{v}+(p-2)\Delta_\infty\hat{v} \right)\leq 1\quad\hbox{ in }x_p.
\end{equation}
Notice that $\nabla v=\nabla\hat{v}$ and $D^2v=D^2\hat{v}$ so  that we can delete the $\hat{} $ s in the last equation. Moreover, the left hand side is continuous in $x_p$. If $|\nabla v(x_\infty)|>1$, then $|\nabla v(x_p)|^{p-4}\to \infty$, while the other factor on the left remains bounded, so that necessarily
$$|\nabla v(x_p)|^2\Delta v(x_p)+(p-2)\Delta_\infty v(x_p)\geq 0$$
as $p\to \infty$ or $-\Delta_\infty v(x_\infty)\leq 0$. The other possibility is $|\nabla v(x_\infty)|\leq 1$. In any case $\min\{ |\nabla v(x_\infty)|-1, -\Delta_\infty v(x_\infty)\leq 0$. This proves that $u_\infty$ is a viscosity subsolution of (\ref{torsioninfty}). 

To see that it is also a viscosity supersolution, suppose that a $C^2$ function $w$ touches $u_\infty$ from below in $x_\infty\in\Omega$. Then by analogous construction $\hat{w}$ touches $u_p$ from below in $x_p$ and $x_p\to x_\infty$ and 
\begin{equation}
-|\nabla w|^{p-4}\left[|\nabla w|^2\Delta w+(p-2)\Delta_\infty w \right]\geq 1\quad\hbox{ in }x_p.
\end{equation}
If $|\nabla w(x_\infty)|<1$, then the left hand side of this relation would tend to zero as $p\to \infty$, a contradition. Therefore $|\nabla w(x_\infty)\geq 1$. Moreover $-[\ldots]\geq 0$ as $p\to\infty$ and thus $-\Delta_\infty w(x_\infty)\geq 0$. but then also $\min\{ |\nabla v(x_\infty)|-1, -\Delta_\infty v(x_\infty)\}\geq 0$. This completes the proof that $u_\infty$ is a viscosity solution of (\ref{torsioninfty}).

One should notice that this proof also provides the fact that $u_\infty$ is a viscosity supersolution to $|\nabla u|=1$. The verification that it is also a viscosity subsolution of $|\nabla u|=1$ requires extra efforts because then the case $|\nabla v(x_\infty)|>1$ must be ruled out, see \cite{BBM}. In \cite{K0} the function $u_\infty$ is identified as $u_\infty(x)=d(x,\partial \Omega)$. We should also mention that by Theorem 2.1 in \cite{Je} solutions to problem (\ref{torsioninfty}) are unique. As it happens, $d(x,\partial\Omega)$ is also a (unique) viscosity solution of the eikonal equation $|\nabla u|=1$ in $\Omega$, $u=0$ on $\partial\Omega$. This observation provides a proof different from \cite{BBM}.

Incidentally, explicit solutions can also be given for the problem $-\Delta_\infty u=1$ in $B_R(0)$. $u=0$ on $\partial B_R(0)$. Since there is uniqueness, solutions are radial and then  in polar coordinates the problem reads $u_r^2u_{rr}=-1$, $u_r(0)=0$, $u(R)=0$ and has the solution $u(r)=c(R^{4/3}-r^{4/3})$ with $c=3^{4/3}/4$. Like the $\infty$-harmonic function $x^{4/3}-y^{4/3}$ this solution is of class $C^{1,1/4}$, and this is the conjectured optimal regularity for equations of the type $-\Delta_\infty u=f$. 

The limit $p\to 1$ in (\ref{ptorsion}) leads to an interesting geometric problem. Formally one looks for minimizers of $J_1(u)=\int_\Omega (|\nabla u|-|u|)\ dx$ and since $u$ can be chosen nonnegative, by the coarea formula and Cavalieri's principle we integrate perimeter minus volume of level sets over all levels:
$$J_1(u)=\int_0^\infty\left( |\partial\{u(x)>t\}| -|\{u(x)>t\}|\right)\  dt.$$
Since level sets are subsets of $\Omega$, a natural question arises. Which subset $\chi_\Omega$ minimizes $|\partial D|-|D|$ or equivalently $\tfrac{|\partial D|}{|D|}$ among subsets of $\Omega$? This is a geometric variational problem, and it has solutions which I call Cheeger sets of $\Omega$. In fact, $h(\Omega):=|\partial \chi_\Omega|/|\chi_\Omega|$ is known as Cheeger constant, and the following estimate of the first Dirichlet Laplace eigenvalue $\lambda$ from below is due to Cheeger: $4\lambda(\Omega)\geq (h(\Omega))^2$, see \cite{Ch,KF,P2}. For convex $\Omega$ Cheeger sets are unique, for nonconvex $\Omega$ they are in general not unique. For plane convex polygons they can be  constructed as in \cite{KL} by sweeping $\Omega$ with discs of suitable radius.  The Cheeger set of a square is explicitly known \cite{K0,KLa}, it is a square with rounded corners of curvature $1/h$; but for the one for a cube in three dimensions there are only numerical approximations \cite{LO}. It is a rounded cube, where the rounded part of the boundary has mean curvature $1/h$.


\section{Parabolic equations}\setcounter{equation}{0}
A parabolic counterpart to $-\Delta_p u=0$ is $u_t-\Delta_pu=0$. This and other equations have been used in mathematical image processing, where typically a blurred black and white picture is represented by a gray scale function $u_0(x)$ defined on a rectangle $\Omega$. Through a nonlinear diffusion filter the function is supposed to evolve into a clearer picture. A popular equation is the so-called total-variation or TV-flow (in which $p=1$). Formally classical solutions of $u_t-{\rm div}(\tfrac{\nabla u}{|\nabla u|})=0$ in $\Omega\times(0,\infty)$ with initial datum $u_0$ and under Neumann boundary conditions have the property that the total spatial variation of $u$ decreases in time. This reduces noise in the picture. But other equations have also been suggested.

In \cite{K1} I made the  Ansatz $v(t,x)=T(t)u(x)$ in $u_t-Au$ for a general operator  $A$ that is homogeneous of degree $d$. This leads to eigenvalue value problems $Au=\lambda u$
and a time decay $T(t)$, which is exponential iff the eigenvalue problem is homogeneous  of degree 1, as in $$Au=-\Delta_p^Nu\qquad\hbox{ or}\qquad Au=-|u|^{2-p}\Delta_pu.$$ 

This is how one arrives at the eigenvalue problems
$$-\Delta_p^Nu=\lambda u\quad\hbox{and}\quad-\Delta_pu=\lambda |u|^{p-2}u \qquad\hbox{ for }p\in(1,\infty),$$
which are treated in subsequent sections. The parabolic equation 
\begin{equation}\label{parabolicnlap}
u_t-\Delta_p^Nu=0
\end{equation}
has interesting special cases. For $p=1$ it is the level set formulation of mean curvature flow and it was extensively studied in a series of papers by Evans and Spruck starting with \cite{ES}. For $p=\infty$ see \cite{JK}. For general $p\in(1,\infty)$ equation (\ref{parabolicnlap}) was studied in \cite{D, BG, JS}.


\section{Dirichlet eigenvalues for $-\Delta_p$}\setcounter{equation}{0}
Consider the Dirichlet eigenvalue problem for the $p$-Laplacian
\begin{equation}\label{pevp}
-\Delta_p u=\lambda_p^p |u|^{p-2}u \quad \hbox{in }\Omega\qquad u=0\quad \hbox{on }\partial\Omega
\end{equation}
with (first nontrivial eigenvalue) $\lambda_p^p$ minimizing 
\begin{equation}
{\cal R}_p(v):=\frac{\int_\Omega |\nabla v|^p}{\int_\Omega |v|^p}\mbox{ on }W^{1,p}_0(\Omega).
\end{equation} 
What is known about $\lambda_p$ and what happens as $p\to\infty$ ?
$\lambda_p\to\lambda_\infty:=1/R(\Omega),$ where $R(\Omega)$ is the inradius of $\Omega$, that is the radius of the largest ball in $\Omega$. 

Moreover, up to a subsequence, the (nonnegative) eigenfunctions $v_p$ converge uniformly to a positive Lipschitzfunction $v_\infty$ solving
\begin{equation}\label{inftyevp}
\min\{|\nabla v|-\lambda_\infty v, -\Delta_\infty v\}=0 \ \mbox{ in }\Omega, \quad v=0 \ \mbox{ on }\partial\Omega
\end{equation}
in the sense of viscosity solutions \cite{JLM}.
 While for finite $p$ the eigenfunctions $v_p$ are unique modulo scaling \cite{KL}, this is no longer the case for $p=\infty$, see \cite{Y,HSY} 

\begin{Remark}\label{rem41}
It is instructive to look at the one-dimensional case in (\ref{inftyevp}), i.e.
\begin{equation}\label{1dinftyevp}
\min\{|v'|-\lambda_\infty v, -|v'|^2v''\}=0 \mbox{ in }\Omega=(-1,1), \quad v(\pm 1)=0 .
\end{equation} 
For $\Omega=(-1,1)$ we have $R(\Omega)=1$, $\lambda_\infty=1$ and $v=1-|x|$ as eigenfunction. Let us see why.

For $|x|\in (0,1)$ any $C^2$-testfunction $\phi$ touching $v$ from above in $x$ satisfies $|\phi'(x)|-\phi(x)=1-v(x)=|x|>0$ and $-\phi''(x)\leq 0$; and in $x=0$ it might have $-\phi''(0)>0$, but in any case $|\phi'(0)|\leq1=\phi(0)$, so $v$ is a viscosity subsolution.

In a similar way one can look at $C^2$-functions $\psi$ touching $v$ from below in $x$. For $x=0$ there are no such functions, so $v$ is also a viscosity supersolution.
\end{Remark}

\begin{figure}[h]
{\vskip4.5cm}
\centerline{\begin{picture}(10,1)(0,0)
\put(3.5,4.5){$\phi$}
\put(5.5,3.25){$\psi$}
\put(8,2){$v$}
\put(0,0){\includegraphics[width=10cm]{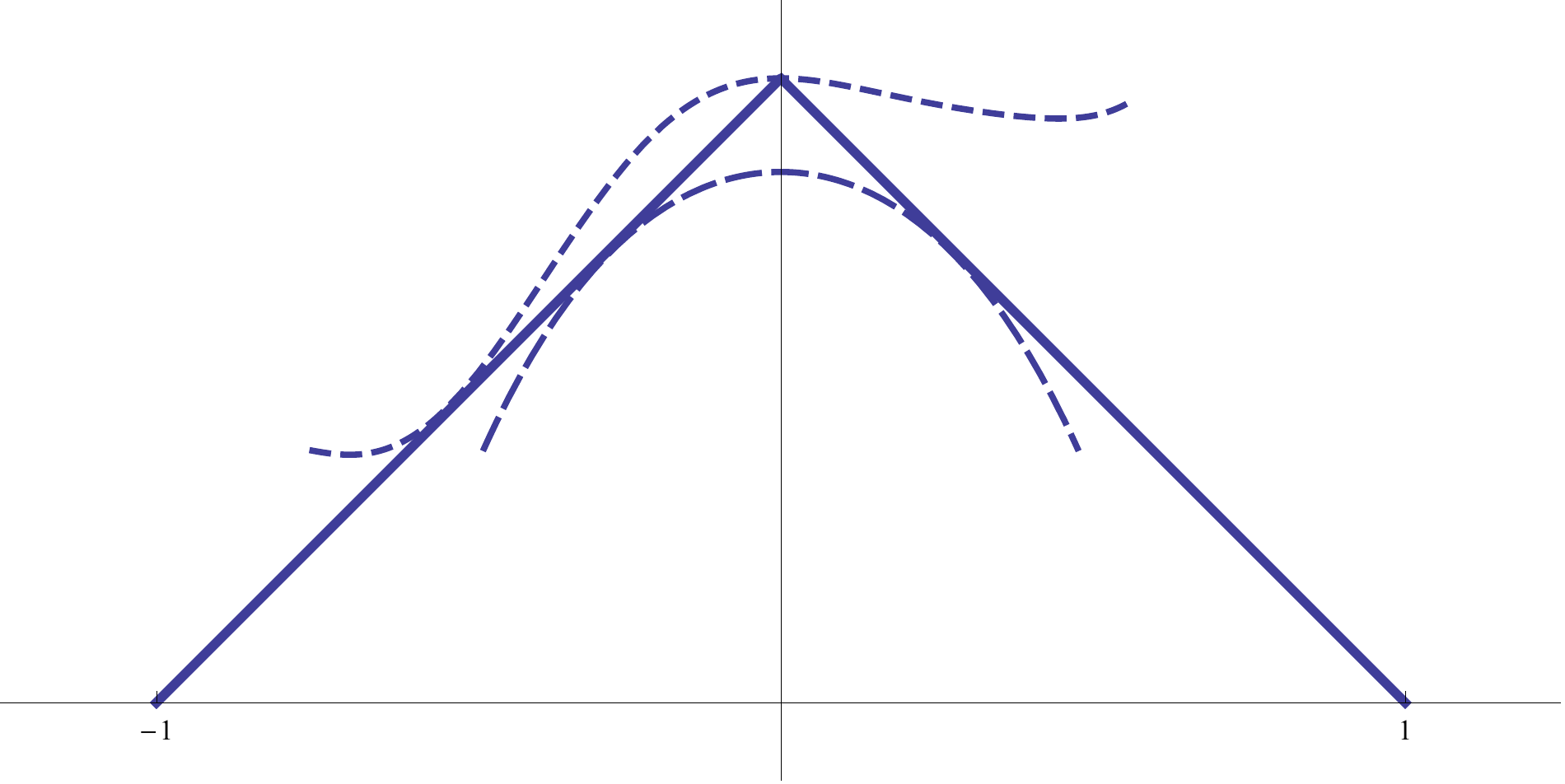}}
\end{picture}}
\caption{The positive viscosity solution of (\ref{1dinftyevp})}
\end{figure}
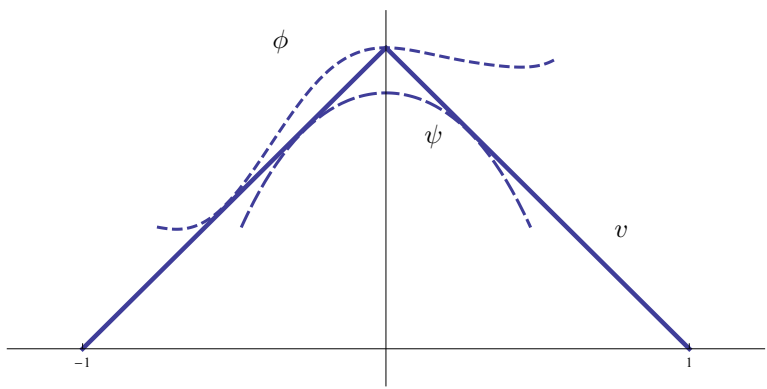

The limit $p\to 1$ in (\ref{pevp}) leads to the formal equation
$$-{\rm div}\left(\frac{\nabla v}{|\nabla v|}\right)=\lambda_1 \left(\frac{ v}{| v|}\right)\quad\hbox{ in }\Omega,
$$
which needs explanation in points where $v$ or $\nabla v$ vanish. If we scale the positive eigenfunctions $v_p$ for $p>1$ in $L^\infty(\Omega)$ to 1, they converge in $BV(\Omega)$ for $p\to 1$ to characteristic functions of a Cheeger set of $\Omega$, see \cite{KF}. Since Cheeger sets are not unique, uniqueness of the limit $v_1$ is probably violated for nonconvex $\Omega$. In \cite{KSc} it was shown that these characteristic functions are indeed viscosity solutions of the limiting equation, but also that there are many more viscosity solutions. 

\begin{figure}[h]
{\vskip3cm}
\centerline{\begin{picture}(10,1)(0,0)
\put(0,0)
{\includegraphics[width=4cm]{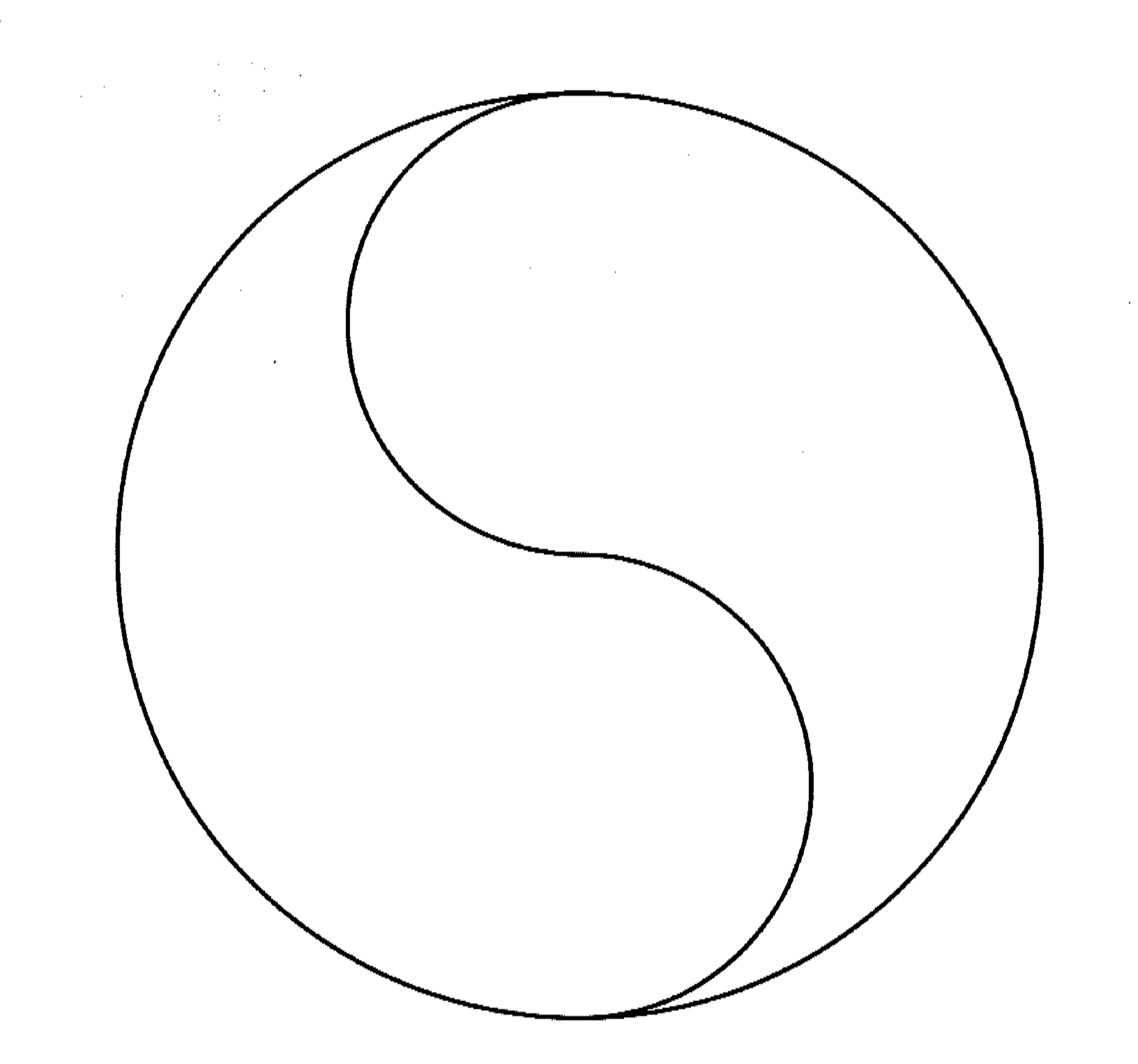}}\hskip6truecm {\includegraphics[width=4cm]{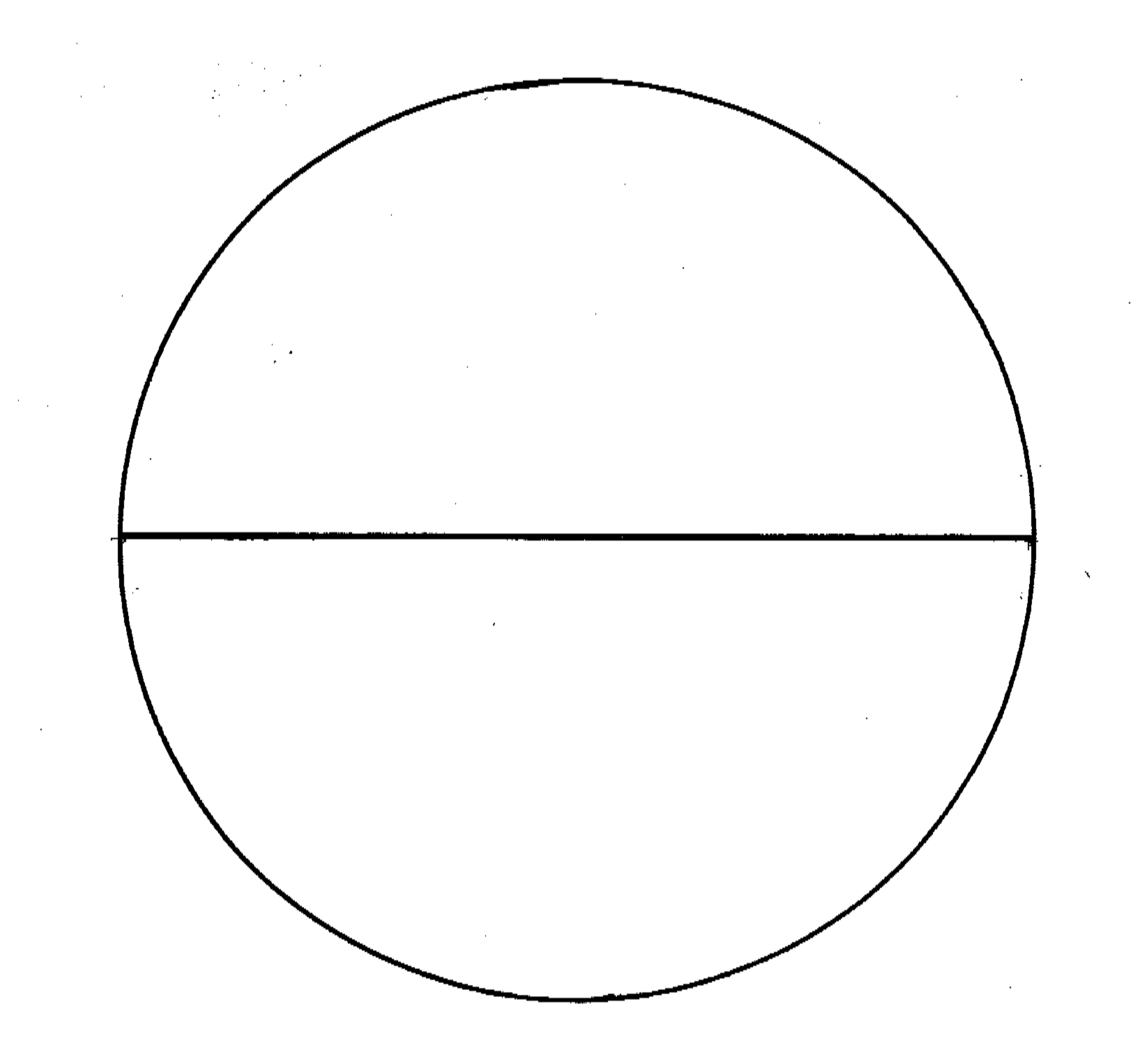}}
\end{picture}}
\caption{Conceivable nodal lines of the second eigenfunction for $p=\infty$ in the disc}
\label{yinyang}
\end{figure}

Results for higher eigenvalues are scarce, partly because beyond the second one there are different conceivable ways of defining them. For $p=\infty$ I refer to \cite{JL}, for $p\to 1$ one can find results in \cite{P}, for $p\in(1,\infty)$ there is a nice numerical study in \cite{H}. Since only the first eigenfunction has one sign, already the second one has a nodal line. Even in case of the unit disc $\Omega=B_1(0)\subset\mathbb{R}^2$ the obvious conjecture that the nodal line should be a diameter is wide open \cite{H,ADS}. For $p=\infty$ it might also be reminiscent of a YinYang-type shape depicted in Figure \ref{yinyang}. 

For a unit square and $p\in(1,2]$ the second eigenfunction seems to have a nodal line parallel to the sides, but for $p\in[2,\infty]$ it appears to be diagonal, see \cite{H,K}.


\section{Neumann eigenvalues for $-\Delta_p$}\setcounter{equation}{0}

Let us now consider the Neumann eigenvalue problem for the $p$-Laplacian
\begin{equation}\label{npevp}
-\Delta_p u=\Lambda_p^p |u|^{p-2}u \quad \mbox{ in }\Omega,\qquad|\nabla u|^{p-2}\frac{\partial u}{\partial n}=0 \quad \mbox{ on }\partial\Omega,
\end{equation}
with (first nontrivial eigenvalue) $\Lambda_p^p$ minimizing 
\begin{equation}
{\cal R}_p(v):=\frac{\int_\Omega |\nabla v|^p}{\int_\Omega |v|^p}\mbox{ on }W^{1,p}(\Omega)\cap\left\{\int_\Omega|v|^{p-2}v =0\right\}.
\end{equation} 
What is known about $\Lambda_p$ and what happens as $p\to 1$ or $p\to\infty$ ?

As $p\to 1$, formally we minimize $\int_\Omega|Du|dx$ on functions satisfying $||u||_{L^1}=1$ and $\int_\Omega {\rm sign} u\ dx=0$.
Minimizers can be chosen constant $\pm c>0$ on disjoint subsets $\Omega^+$ and $\Omega^-$ of $\Omega.$

This constitutes a geometric partioning problem: 
Divide $\Omega$ into two disjoint subsets $\Omega^+$ and $\Omega^-$ of equal volume 
such that their relative perimeter $P(\partial\Omega^+;\Omega)$ becomes minimal.
Problems of this nature are relevant in developing fast algorithms 
for parallel computing \cite{GG}.


For $p\to\infty$ we have $\Lambda_p \to \Lambda_\infty:=2/diam(\Omega)$, while $\lambda_\infty=1/R(\Omega)$.
The first nontrivial eigenfunctions $u_p$ converge to a viscosity solution $u_\infty$ of \begin{equation}\label{neui}\begin{cases}
 {\min}  \{| \nabla u| - \Lambda u, -\Delta_{\infty} u\} = 0 & \mbox{in}\ \{u >0\}\cap\Omega\\
{\max}  \{-| \nabla u|  - \Lambda u , -\Delta_{\infty} u\} =0 & \mbox{in}\ \{u <0\}\cap\Omega\\
-\Delta_{\infty} u  =0 & \mbox{in} \ \{u =0\}\cap\Omega\\
{\partial u}/{\partial \nu} =0 & \mbox{on} \ \partial\Omega.\\
\end{cases}
\end{equation}
While (\ref{neui}) was already derived in \cite{RS}, we also prove in \cite{EKNT} that $\Lambda\geq\Lambda_\infty$ for any $\Lambda$ having a nontrivial solution to (\ref{neui}) and any convex $\Omega$. This shows that (for convex $\Omega$) $\Lambda_\infty$ is in fact the first positive eigenvalue of (\ref{neui}).

The fact that $\Lambda_\infty=2/diam(\Omega)$ has interesting consequences.
\begin{cor}\label{SW}
 If $\Omega^*$ denotes a ball of same volume as $\Omega$, then the Szeg\"o-Weinberger inequality $\Lambda_\infty(\Omega)\leq \Lambda_\infty(\Omega^*)$ holds.
\end{cor}
For general $p$ we do not know that $\Lambda_p(\Omega)\leq \Lambda_p(\Omega^*)$, except

a) for $p=2$ and general $\Omega$ or 

b) for $p=\infty$ and convex $\Omega$ (Corollary \ref{SW}).

c) For $p=1$ and convex plane $\Omega$ this was an old conjecture of Poly\'a, which was finally solved after 60 years in \cite{EFKNT}. The proof is quite involved.  
\begin{cor}\label{DN}
For convex $\Omega$ we have $\Lambda_\infty(\Omega)\leq\lambda_\infty(\Omega)$.  Moreover,  equality holds only if $\Omega$ is a ball.
\end{cor}
A stronger statement is known for finite $p$, namely $\Lambda_p(\Omega)<\lambda_p(\Omega)$

a) for $p=2$ and general $\Omega$ as a consequence of the Faber-Krahn and Szeg\"o-Weinberger inequalities.

b) for $p\in (1,\infty)$ and convex $\Omega$ \cite{BNT}.
\begin{cor}\label{closednodal} For convex $\Omega$ no Neumann eigenfunction associated to $\Lambda_\infty$ 
 can have a closed nodal line inside $\Omega$.
\end{cor}
In fact, otherwise there is a nodal domain $\Omega '\subset\subset\Omega$ so that $\Lambda_\infty(\Omega)=\lambda_\infty(\Omega ')=1/R(\Omega ')>1/R(\Omega)=\lambda_\infty(\Omega)$ contradicting Corollary {DN}

Let us now see that $\Lambda_p\to\Lambda_\infty$, and where convexity of $\Omega$ enters into the proof that $\Lambda_\infty$ is the first nontrivial Neumann eigenvalue in the limiting problem (\ref{neui}). 

\begin{lemma}\label{lemma1}
For $\Omega$ simply connected and Lipschitz $\Lambda_p\to\Lambda_\infty:=2/diam(\Omega)$, where $diam$ demotes intrinsic diameter.
\end{lemma}
\smallskip

For $x, y\in\Omega$ the intrinsic distance $d_\Omega(x,y)$ is the length of the geodesic in $\Omega$ connecting $x$ to $y$ and $diam(\Omega)=\sup_{x,y\in\Omega}d_\Omega(x,y)$.

The proof of Lemma \ref{lemma1} is done in two steps.

\noindent {\bf Step 1:} $\limsup_{p\to\infty}\Lambda_p\leq\Lambda_\infty$.
Pick $x_0\in\Omega$ and adjust $c_p\in\mathbb{R}$ so that $w(x)=d_\Omega(x,x_0)-c_p$ is admissible test function for ${\cal R}_p$. Then
$$\Lambda_p\leq{\cal R}_p(w)=\left(\tfrac{1}{|\Omega|}\int_\Omega |d_\Omega(x,x_0)-c_p|^p\right)^{-1/p}.$$
Since $0\leq c_p\leq diam(\Omega)$, up to a subseqence $c_p\to c_\infty$ and
$$\liminf_{p\to\infty}\left(\tfrac{1}{|\Omega|}\int_\Omega |d_\Omega(x,x_0)-c_p|^p\right)^{1/p}=sup_{x\in\Omega} |d_\Omega(x,x_0)-c_\infty|\geq diam(\Omega)/2=\frac{1}{\Lambda_\infty}.$$

\noindent {\bf Step 2:} $\liminf_{p\to\infty}\Lambda_p\geq\Lambda_\infty$.
For $p>m>n$ the family of eigenfunctions $u_p$ is uniformly bounded in $W^{1,m}(\Omega)$ and equicontinuous so it converges (after possibly passing to a subsequence) in $C^0(\Omega)$ and weakly in $W^{1,m}(\Omega)$ to a limit $u_\infty$. Thus
$$\frac{||\nabla u_\infty||_m}{||u_\infty||_m}\leq\liminf_{p\to\infty}\frac{||\nabla u_p||_m}{||u_p||_m}
\leq
\liminf_{p\to\infty}\frac{||\nabla u_p||_p}{||u_p||_m}
=\liminf_{p\to\infty}\Lambda_p\frac{|| u_p||_p}{||u_p||_m}.$$
Sending $m\to \infty$ gives
$\liminf_{p\to\infty}\Lambda_p\geq\frac{||\nabla u_\infty||_\infty}{||u_\infty||_\infty}$, and  
so it remains to estimate $||u_\infty||_\infty$ in terms of $||\nabla u_\infty||_\infty$.

But condition $\int_\Omega |u_p|^{p-2}u_p=0$ implies $\sup u_\infty=-\inf u_\infty$
 $$||u_\infty||_\infty=\tfrac{1}{2}(\sup u_\infty-\inf u_\infty)\leq \tfrac{1}{2} diam(\Omega) ||\nabla u_\infty||_\infty .$$
Thus $$\liminf_{p\to\infty}\Lambda_p\geq\frac{||\nabla u_\infty||_\infty}{||u_\infty||_\infty}\geq\frac{2}{diam(\Omega)}=\Lambda_\infty$$

\hfill\qed

\begin{Remark} Suppose $||u_\infty||_\infty$ is scaled to 1, $u(x)=\inf u$ and $u(y)=\sup u$.
Then $u_\infty$ increases with constant slope $\Lambda_\infty$ along the geodesic from $x$ to $y$.

In particular in the one-dimensional case that $\Omega=(-1,1)$ the function $u(x)=x$ is a viscosity solution of (\ref{neui}). While it solves the differential equation $-\Delta_\infty u=u''=0$ in the classical (and thus in the viscosity) sense, it clearly violates the boundary condition $u'(\pm1)=0$ in the classical sense. To check the Neumann condition in the right end point $x=1$ point in the viscosity sense {\rm \cite{CIL}}, one must 
 must verify
\begin{equation}\label{nc1}
\min \{ \min \{|\phi'|-\Lambda \phi, -|\phi'|^2\phi''\}\ , \phi'\}(1)\leq 0
\end{equation}
for any $C^2$ test function $\phi$ touching $u$ in $x=1$ from above, and
\begin{equation}\label{nc2}
\max\{ \min \{|\psi'|-\Lambda \psi, -|\psi'|^2\psi''\}\ , \ \psi'\}(1)\geq 0
\end{equation}
for any smooth test function $\psi$ touching $u$ from below.

\begin{figure}[h]
{\vskip8cm}
\centerline{\begin{picture}(10,0)(0,0)
\put(6.5,7.3){$\phi_1$}
\put(6.3,6){$\phi_2$}
\put(7.5,4.2){$\psi$}
\put(5.5,4.2){$v$}
\put(1.5,0){\includegraphics[width=7cm]{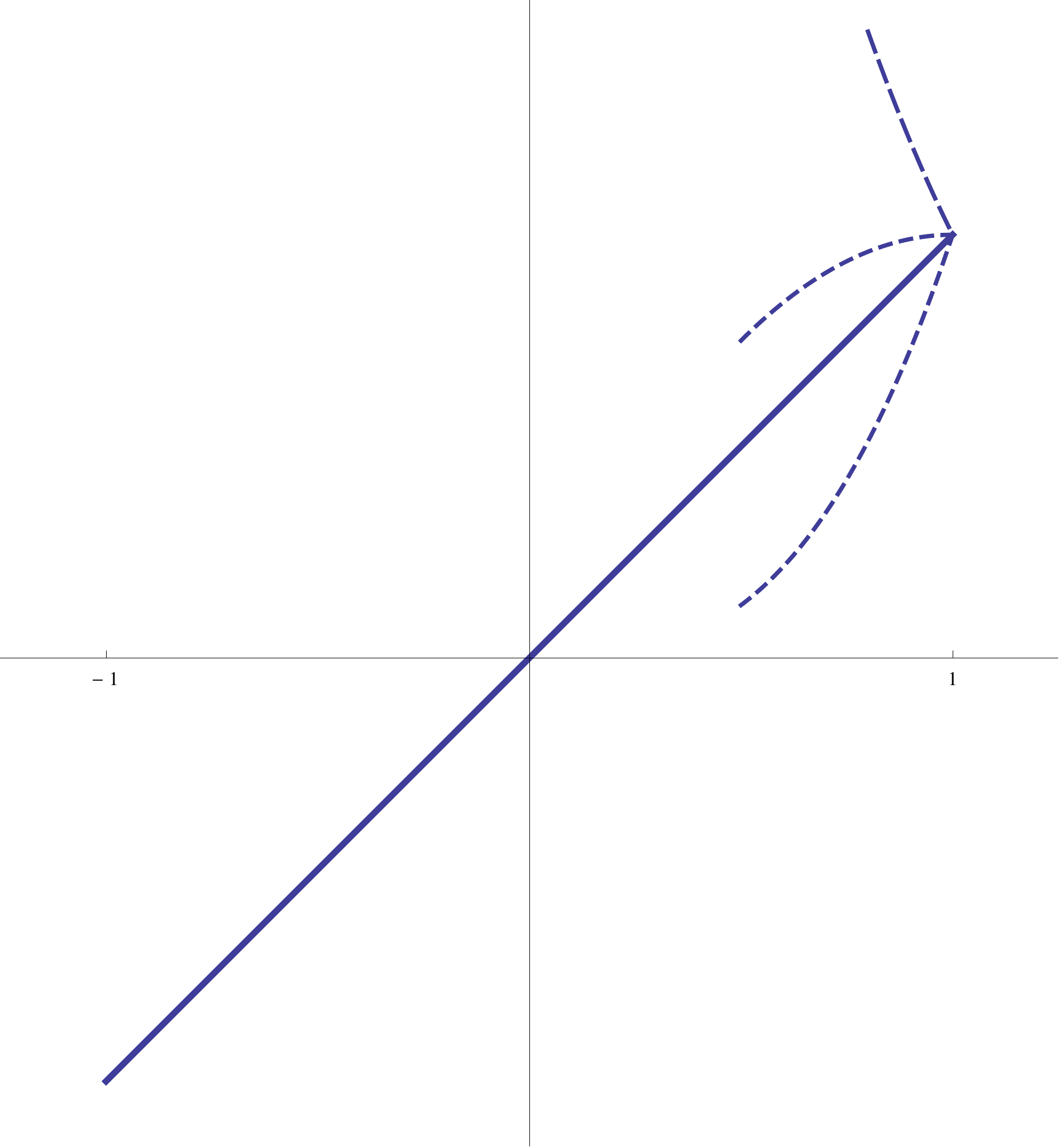}}
\end{picture}}
\caption{Illustration of (\ref{nc1}) and (\ref{nc2})}
\end{figure}
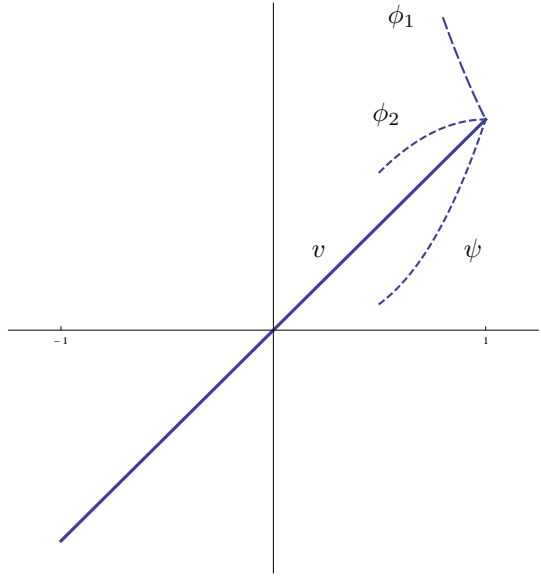 

\end{Remark}

The fact that boundary conditions are not satisfied in a classical sense is quite common for viscosity solutions. In \cite{KK}, for instance, it was pointed out that the problem $u_t-\Delta u=0$ in $\Omega\times (0,T)$ has a viscosity solution satisfying initial data $u(x,0)=u_0$ in $\Omega$, boundary data $u(x,t)=g(x,t)$ on $\partial\Omega\times(0,T)$ {\bf and} final data $u(x,T)=u_T(x)$ in $\Omega$ for arbitrary continuous $u_T$ compatibel with the boundary datum $g$. The final condition is in general violated in the classical sense, but it holds in the viscosity sense, because the solution of the initial boundary value problem is classical at time $T$.

\begin{figure}[h]
  \centerline{\includegraphics[width=5cm,clip=true,viewport= 300 200 2000 1600]{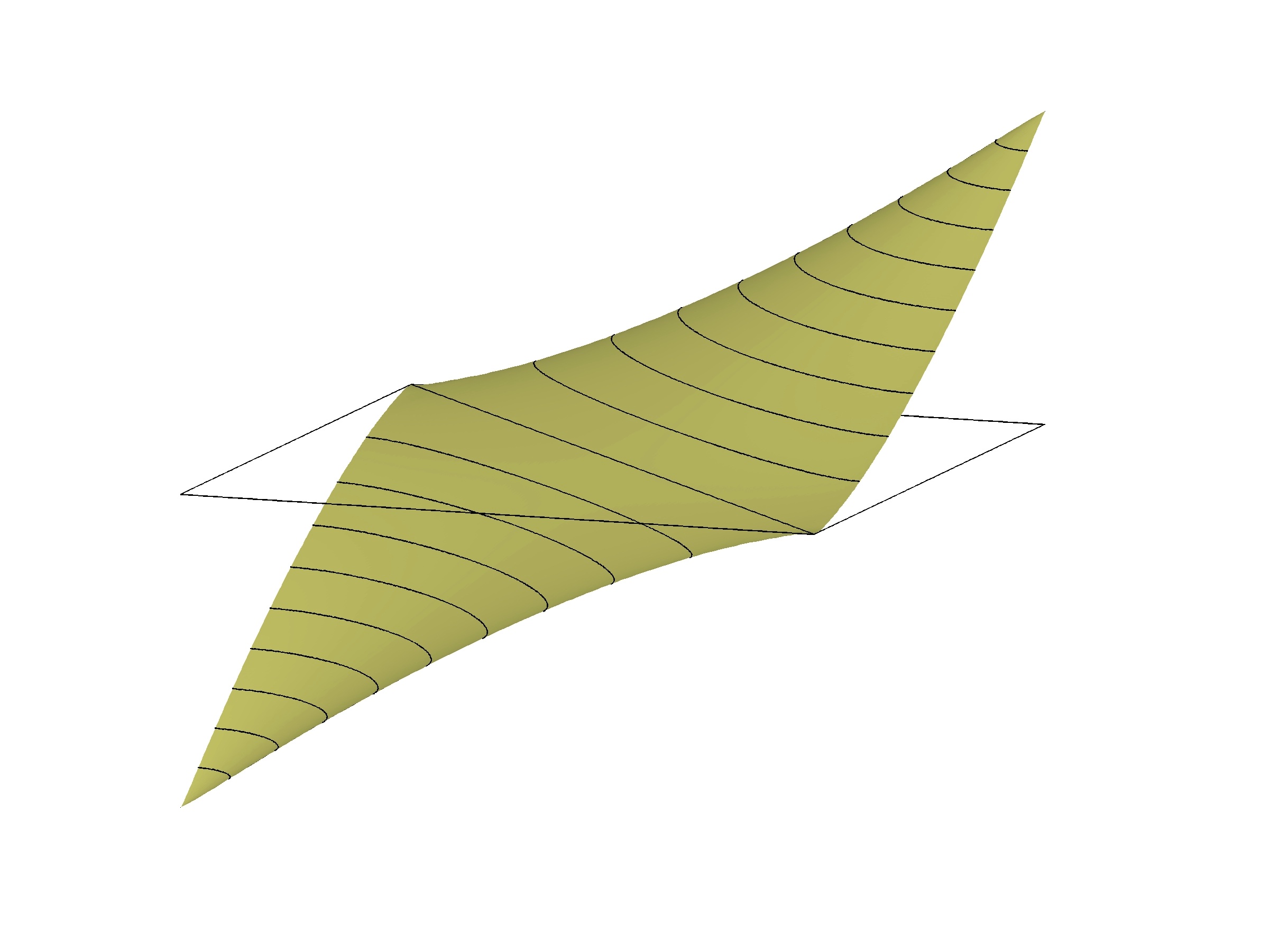}\hspace{10mm}
    \includegraphics[width=5cm,clip=true,viewport= 300 150 2200 1600]{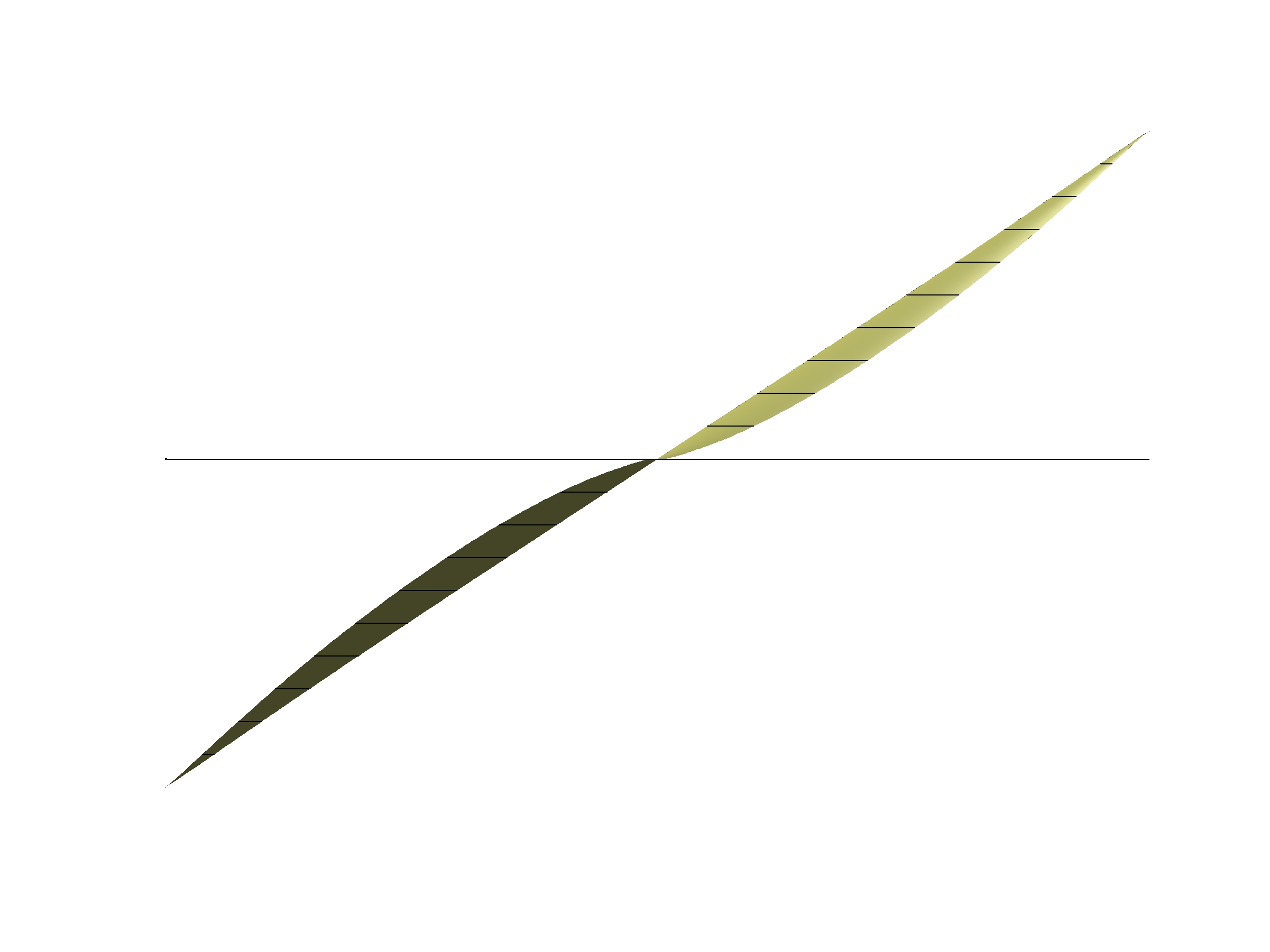}}
  \caption{Numerical simulation of $u_{15}$ and side view in diagonal direction.}
  \label{Neumann_infty_on_square1}
\end{figure}

\begin{figure}[h]
  \centering
  \setlength{\unitlength}{1mm}
  \begin{picture}(120,45)
    \put(0,0){\includegraphics[width=60mm]{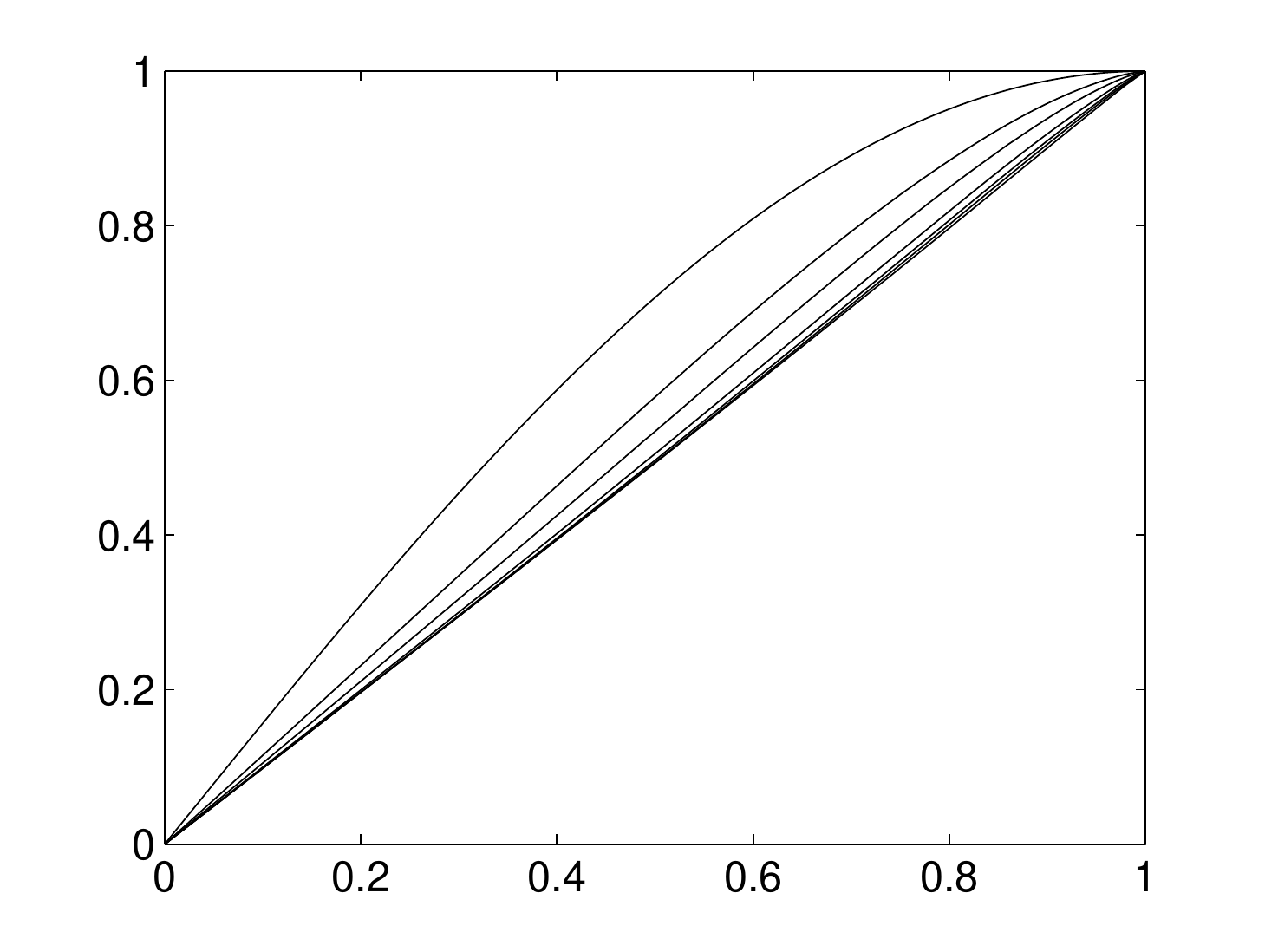}}
    \put(60,0){\includegraphics[width=60mm]{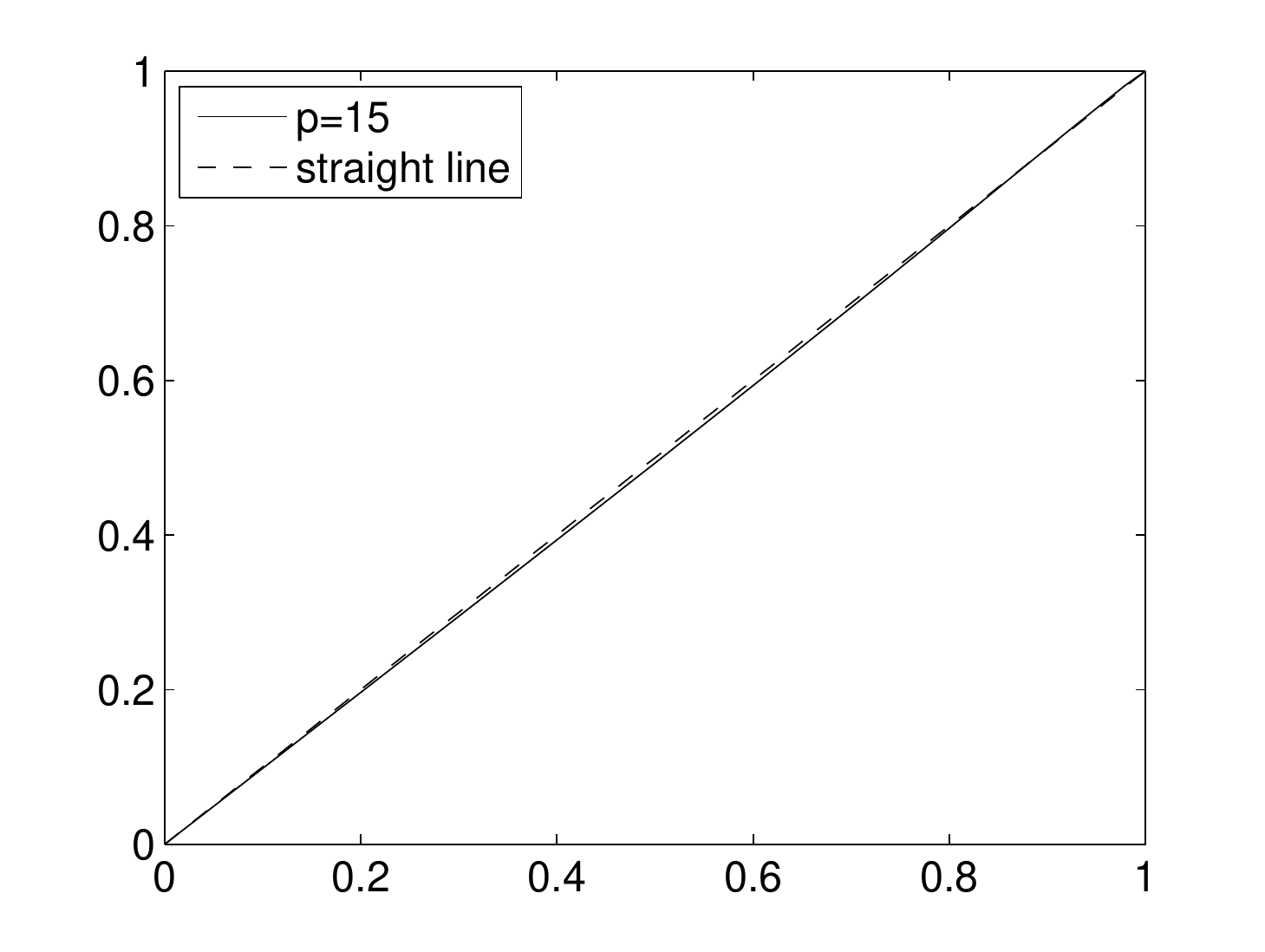}}
    \put(101,8){\includegraphics[width=20mm,clip=true,viewport= 88 30 350 300]{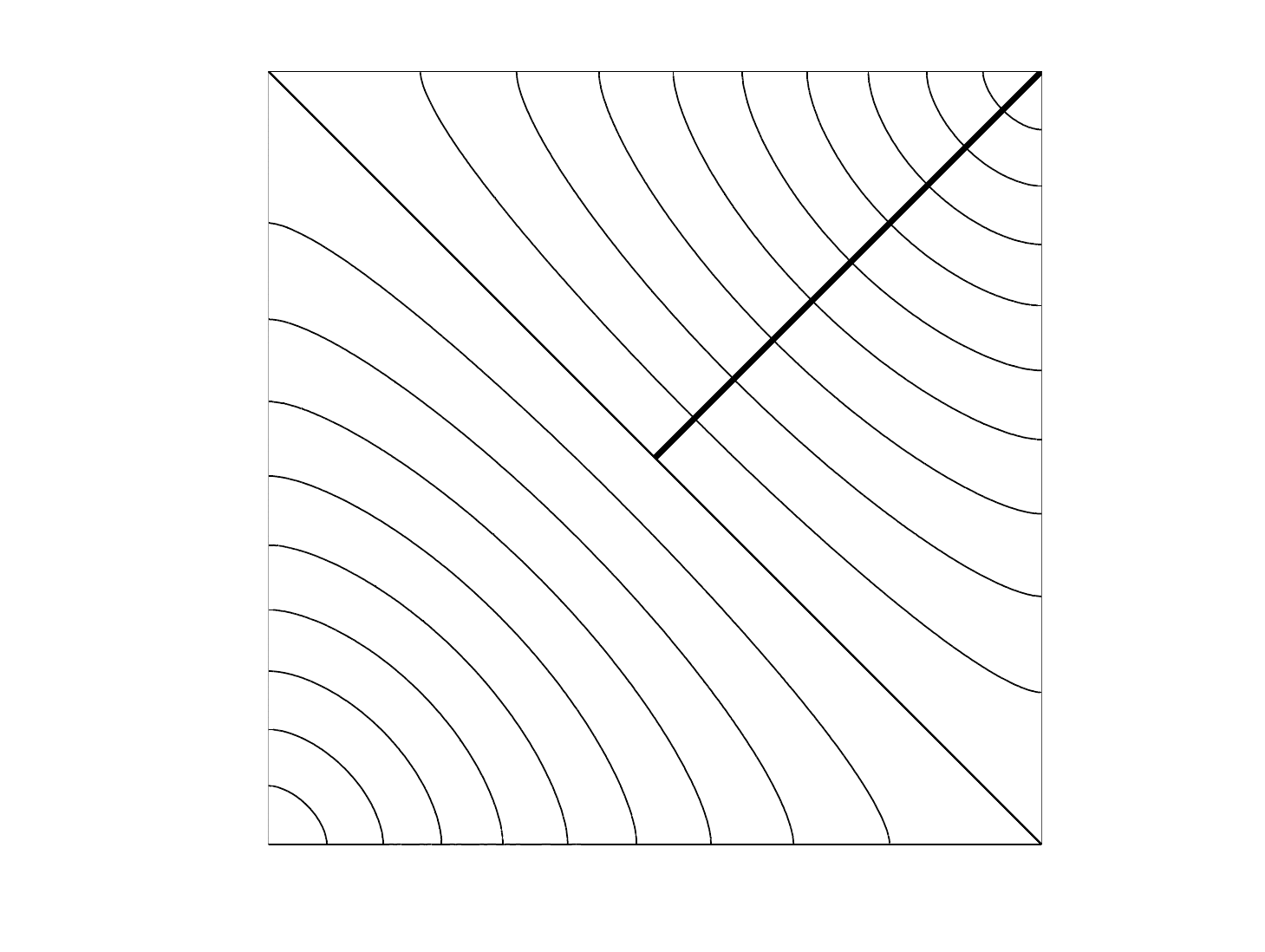}}
  \end{picture}
  \caption{Numerical simulation of $u_p$: normalized values along half
    of the diagonal for $p=2, 3, 4, 6, 8, 10, 15$ (left), and for $p=15$
    compared to the line $y=x$ (right).}
  \label{Neumann_infty_on_square2}
\end{figure}

By the previous remark, on a rectangle $u_\infty$ has constant (and maximal) slope along a diagonal, and numerical simulations of the first nonconstant eigenfunction $u_p$ for increasing $p$ as in Figures \ref{Neumann_infty_on_square1} and \ref{Neumann_infty_on_square2} confirm this behaviour. The computations were done by adapting the steepest descent method of \cite{H}.

\begin{Remark}
For a square and the classical Laplacian $\Lambda_2$ is a double eigenvalue and there are solutions with diagonal and horizontal or vertical nodal lines, and for a rectangle that is not a square and $p=2$, the associated (linear) eigenfunction $u_2$ has its gradient parallel to the longer sides. 
It is instructive to note that there are one-dimensional solutions also for {\rm (\ref{neui})} on rectangles; but they are associated to higher eigenvalues than $\Lambda_\infty$. 
Take $\Omega=(-1,1)^2$ and $u(x)=x_1$. Then $u\in C^2(\Omega)$, $-\Delta_\infty u=0$ in $\Omega$. 

Now the PDE in (\ref{neui}) is satisfied  if also $1=|\nabla u|\geq\Lambda u$ on $\{ u>0\}$, that implies $$\Lambda\leq 1.$$

The Neumann boundary condition is satisfied in the classical sense on horizontal parts of  $\partial\Omega$ and violated in the classical sense on vertical parts of $\partial\Omega$. 

However, for the Neumann boundary condition to hold in the viscosity sense (see {\rm \cite{CIL}}) on the right part,  we must verify
$$\min \{ \min \{| \nabla \phi|-\Lambda \phi, -\Delta_\infty\phi\}\ , \ \partial\phi/\partial \nu\}(x_0)\leq 0$$
for any $C^2$ test function $\phi$ touching $u$ in $x_0\in\partial\Omega$ from above, and
$$\max\{ \min \{|\nabla \psi|-\Lambda \psi, -\Delta_\infty\psi\}\ , \ \partial\psi/\partial \nu\}(x_0)\geq 0$$
for any smooth test function $\psi$ touching $u$ from below. 

Recall $|\nabla u|=\partial u/\partial\nu=1$ here. Therefore only the very first constraint is active and implies $$\Lambda\geq 1.$$
This shows that $u(x)=x_1$ is a viscosity solution to {\rm (\ref{neui})} with eigenvalue $\Lambda=1$, but $$\Lambda=1>\frac{1}{\sqrt{2}}=\frac{2}{diam(\Omega)}=\Lambda_\infty.$$
\end{Remark}

So far we have shown  that $\Lambda_p \to \Lambda_\infty:=2/diam(\Omega)$ and by standard arguments as in \cite{JLM}
one can also see that the first nontrivial eigenfunctions $u_p$  of (\ref{npevp}) converge to a viscosity solution $u_\infty$ of (\ref{neui}). 

Let us now prove that $\Lambda\geq\Lambda_\infty$ 
for any $\Lambda$ having a nontrivial solution to (\ref{neui}) and any convex $\Omega$. 
This shows that (for convex $\Omega$) $\Lambda_\infty$ is in fact the first positive eigenvalue of (\ref{neui}) and $u_\infty$ is at least one associated eigenfunction.

\begin{thm}\label{theo57}
If $u$ is a nontrivial solution to (\ref{neui}) with $\Lambda>0$  
and if $\Omega$ is convex and smooth, then $\Lambda\geq \Lambda_\infty$.
\end{thm}

The proof requires a few intermediate steps.

\begin{lemma} {\rm (Strong maximum principle)} If $\Omega_1\subset\Omega$ is open and connected 
and $u\geq m>0$ on $\overline{\Omega_1}$, then $u>m$ in $\Omega_1$.
\end{lemma} 
In fact $\min\{ |\nabla u|-\Lambda u,-\Delta_\infty\}=0$ in $\Omega_1$, so $-\Delta_\infty u\geq 0$ and $|\nabla u|\geq \Lambda u$ in $\Omega_1$. In particular, if $u\in C^1$ attains a local min inside $\Omega_1$, then $|\nabla u|$ vanishes there, contradicting $|\nabla u|\geq \Lambda m>0$. These heuristics can be replaced by a precise proof.\hfill\qed

\begin{lemma} 
$u$ must change sign.
\end{lemma}

Without loss of generality$u$ is positive somewhere and cannot have positive min in $\Omega$.
If the minimum is zero in $x_0\in\Omega$, then comparison with a flat cone yields a contradiction. Any nonnegative minimum on $\partial\Omega$ contradicts the Neumann boundary condition in the viscosity sense. Therefore the minimum of $u$ must be negative.\hfill\qed

In particular the nodal set $\{ x\in\Omega\ ;\ u(x)=0\}$ is nonempty.
For the proof of Theorem \ref{theo57} we pick a point $x_0$ in $\Omega_-:=\{x\in\Omega\ ; \ u(x)<0\}$,
scale $\max u$ to $1/\Lambda$
and construct a smooth test function (off of $x_0$)
$$g_{\varepsilon,\gamma}:=(1+\varepsilon)|x-x_0|-\gamma|x-x_0|^2.$$ 
It is a viscosity supersolution to (\ref{neui}) in $\Omega_+$ and satisfies $g_{\varepsilon,\gamma}(x)\geq u(x)$ there.
Sending $\varepsilon$ and $\gamma$ to zero and $x_0$ to a nodal point of $u$ gives
$$|x-x_0|\geq u(x)\qquad\mbox{ and }\qquad d^+=\sup_{x\in \Omega_+} dist(x, \{ u=0\})\geq\frac{1}{\Lambda}.$$
Similarly 
$$d^-=\sup_{x\in \Omega_-} dist(x, \{ u=0\})\geq\frac{1}{\Lambda}.$$
$$\hbox{ But }\qquad\qquad diam(\Omega)\geq d^+ +d^-\geq \frac{2}{\Lambda}$$
which completes the proof of Theorem \ref{theo57}. \hfill\qed

For nonconvex $\Omega$ it would be natural to replace $|x-x_0|$ in the construction of $g_{\varepsilon,\gamma}$ by geodesic distance $d_\Omega(x,x_0)$, but unfortunately then $g_{\varepsilon,\gamma}\not\in C^2$ is no longer an admissible testfunction.

\section{Dirichlet problems for $-\Delta_p^Nu=0$ and $-\Delta_p^Nu=1$}
\setcounter{equation}{0}

The Dirichlet problem
\begin{equation}\label{normharm}
-\Delta_p^Nu=0\quad\hbox{ in }\Omega,\qquad u=g\quad\hbox{ on }\partial\Omega
\end{equation}
has been less studied than the version described in Section \ref{Dirichletsection}. While weak solutions to $-\Delta_p u=0$ are viscosity solutions and vice versa \cite{JJ} for any $p\in(1,\infty)$, the relation between viscosity solutions to a) $-\Delta_p u=0$ and b) $-\Delta_p^Nu=0$ appears to be more delicate. As mentioned in \cite{K} for $p\in(1,2)$ any viscosity solution of a) is one of b), and upon inspecting the definitions for $p\in(2,\infty)$ any viscosity solution of b) is also one of a). The situation is similar in the inhomogeneous case
\begin{equation}\label{normtors}
-\Delta_p^Nu=1\quad\hbox{ in }\Omega,\qquad u=0\quad\hbox{ on }\partial\Omega,
\end{equation}
where one can have more solutions to a) $-\Delta_pu=p|\nabla u|^{p-2}$ than to the seemingly equivalent equation
b) $-\Delta_p^Nu=1$. This was exhibited in \cite{CF} in a one-dimensional setting, and for the sake of exposition let us look at it in the case of $\Omega=B_R(0)\subset\mathbb{R}^n$. In this case $u(x)=v(|x|)$ and (\ref{normtors}) turns into
\begin{equation}\label{normtorsode}
-(p-1)v''-\frac{n-1}{r}\ v'=p\quad\hbox{ in }(0,R),\qquad v'(0)=0=v(R),
\end{equation}
with the solution $v(r)=c(n,p)\left(R^2-r^2\right)$ and with $c(p,n)=\frac{p}{2(n-2+p)}$ nondecreasing in $p$. It is curious to note that for $n=2$ this solution is independent of $p$

Anyway, version a) of (\ref{normtorsode}) would read
$$
-(p-1)|v'(r)|^{p-2}v''-\frac{n-1}{r}\ |v'(r)|^{p-2}\ v'=p|v'(r)|^{p-2}\ \hbox{ in }(0,R),\quad v'(0)=0=v(R),
$$
and for any $p>2$ and $\rho\in(0,1)$ it has the family of piecewise $C^2$ solutions
$$v(x):=c(p,n)\begin{cases}
(R-\rho)^2&\hbox{ if }r\in[0,\rho],\cr
(R-\rho)^2-(r-\rho)^2 &\hbox{ if }r\in[\rho,R].
\end{cases}
$$
It is easy to check that these are also viscosity solutions of $-\Delta_pu=p|\nabla u|^{p-2}$, while they cannot be viscosity solutions of $-\Delta_p^Nu=1$ or (\ref{normtorsode}).

Finally let us not forget to mention that $C^1$ regularity of viscosity solutions to the equation $-\Delta_p^N=f$ was derived in \cite{BD,APR, CF}, and that solutions of the corresponding Dirichlet problem are unique for positive $f$ and any $p\in(2,\infty]$ in \cite{BD,LW,KMP}, but only generically unique for $p=\infty$ and $f$ changing sign in \cite{AS} .


\section{Dirichlet eigenvalues for $-\Delta_p^N$}
\setcounter{equation}{0}

The eigenvalue problem
\begin{equation}\label{normevp}
-\Delta_p^Nu=\lambda_p u\quad\hbox{ in } \Omega, \qquad u=0\quad \hbox{ on }\partial\Omega
\end{equation} 
was studied for instance in \cite{BD,MPR} for $p\in(1,\infty)$ and  for $p=\infty$ in \cite{J}. For $p\geq2$ it was shown in \cite{BD2} that the first eigenfunction is again unique modulo scaling. At least for starshaped domains $\Omega$ there is convergence \cite{MPR}  of the first eigenvalue $\lambda_p$ to $\lambda_\infty$ as $p\to\infty$, but the limit $p\to 1$ appears to be open.

Therefore I turn to the special case $\Omega=B_R(0)$ of a ball and to radially symmetric solutions.
One should also expect many nonradial eigenfunctions to exist, but for $p\not=2$ we are not aware of any results in this direction, not even in two dimensions. 
For radial functions  and with the Ansatz $u(x)=v(|x|)$ the eigenvalue problem (\ref{normevp})
 transforms into
\begin{equation}\label{ode} 
(p-1)v''(r)+\frac{n-1}{r}v'(r)+p\lambda v(r)=0 \quad\hbox{ in }(0,R), \qquad
v'(0)=0=v(R).
\end{equation}
The formal limit $p\to 1$ in (\ref{ode}) is
\begin{equation}\label{degode}\frac{n-1}{r}v'(r)+\lambda v(r)=0\quad\hbox{ in }(0,R),\qquad v'(0)=0=v(R).
\end{equation}
 with the viscosity solution $v(r)=v(0)\  e^{-\tfrac{\lambda}{2(n-1)}r^2}$ , which violates the boundary condition in the classical sense at $r=R$, but since it solves (\ref{degode}) there in the classical sense, it still satisfies the boundary condition $v(R)=0$ in the viscosity sense. It is therefore reasonable to expect that also for more general domains the first eigenfunctions in (\ref{normevp}) converge (as $p\to 1$) to a viscosity solution of
 \begin{equation}\label{conjp1}
 (n-1)H v_\nu+\lambda v=0\quad \hbox{ in }\Omega, \qquad v=0\quad\hbox{ on }\partial\Omega,
 \end{equation}
and that they develop a boundary layer.

For $p>1$ problem (\ref{ode}) is a Bessel type equation, and this observation was used in \cite{KKK} to 
explicitly derive a countable and complete orthonormal system of eigenfunctions to (\ref{ode}) in a suitably weighted $L^2$-space. These form a countable, but presumably incomplete sequence of eigenfunctions to (\ref{normevp}).

\bigskip\noindent
Authors' addresses:

\medskip\noindent
Bernd Kawohl
\hfill\break
Mathematisches Institut
\hfill\break
Universit\"at zu K\"oln
\hfill\break
50923 K\"oln, Germany
\hfill\break
kawohl@mi.uni-koeln.de
\medskip

\noindent
Ji\v r\'{\i} Hor\'ak
\hfil\break
Fakult\"at Maschinenbau
\hfill\break
TH Ingolstadt
\hfill\break
Postfach 21 04 54
\hfil\break
85019 Ingolstadt, Germany
\hfill\break
Jiri.Horak@thi.de
\end{document}